\documentclass[11pt,onecolumn]{article}
\usepackage{multicol}
\usepackage[greek,english]{babel}
\usepackage[vcentering,dvips]{geometry}
\usepackage{amsfonts}
\usepackage{amssymb}
\usepackage{eucal}
\usepackage{amsxtra}
\usepackage{color}
\usepackage{setspace}
\usepackage{times}
\usepackage{titlesec}
\usepackage{sectsty}
\allsectionsfont{\large}
\usepackage{url}
\usepackage{doi}
\usepackage{hyperref}
\hypersetup{
    colorlinks,%
    citecolor=black,%
    filecolor=black,%
    linkcolor=black,%
    urlcolor=black
}

\addtocounter{section}{0} \numberwithin{equation}{section}
\newtheorem{theorem}{Theorem}[section]
\newtheorem{proposition}[theorem]{Proposition}
\newtheorem{lemma}[theorem]{Lemma}

\newcommand{\qed}{{$\hfill \Box$}}

\onehalfspace

\hypersetup{urlcolor=blue }

\begin{document}
\pagestyle{plain}

\title{\large
{\textbf{REAL HYPERSURFACES EQUIPPED WITH $\xi$-PARALLEL STRUCTURE JACOBI OPERATOR IN $\mathbb{C}P^{2}$ OR
$\mathbb{C}H^{2}$}}}

\author{ \textbf{\normalsize{Konstantina Panagiotidou and Philippos J. Xenos}}\\
\small \emph{Mathematics Division-School of Technology, Aristotle University of Thessaloniki, Greece}\\
\small \emph{E-mail: kapanagi@gen.auth.gr, fxenos@gen.auth.gr}}
\date{}

\maketitle
\begin{flushleft}
\small {\textsc{Abstract}. The $\xi$-parallelness condition of the structure Jacobi operator of real hypersurfaces has been studied in combination with additional conditions. In the present paper we study three dimensional real hypersurfaces in $\mathbb{C}P^{2}$ or
$\mathbb{C}H^{2}$ equipped with $\xi$-parallel structure Jacobi operator. We prove that they are Hopf hypersurfaces and if additional $\eta(A\xi)\neq0$, we give the classification of them.}
\end{flushleft}
\begin{flushleft}
\small{\emph{Keywords}: Real hypersurface, $\xi$-parallel structure Jacobi operator, Complex projective space, Complex hyperbolic space.\\}
\end{flushleft}
\begin{flushleft}
\end{flushleft}

\section{Introduction}

A complex n-dimensional Kaehler manifold of constant holomorphic sectional curvature c is called a complex space form, which is denoted by $M_{n}(c)$. A complete and simply connected complex space form is complex analytically isometric to a complex projective space $\mathbb{C}P^{n}$, a complex Euclidean space $\mathbb{C}^{n}$ or a complex hyperbolic space $\mathbb{C}H^{n}$ if $c>0, c=0$ or $c<0$ respectively.

The study of real hypersurfaces in a nonflat complex space form is a classical problem in Differential Geometry. Let $M$ be a real
hypersurface in $M_{n}(c)$. Then $M$ has an almost contact metric structure $(\varphi,\xi,\eta,g)$. The structure vector field $\xi$
is called principal if $A\xi=\alpha\xi$ holds on $M$, where A is the shape operator of $M$ in $M_{n}(c)$ and $\alpha$ is a smooth function. A real hypersurface is called \textit{Hopf hypersurface} if $\xi$ is principal.

Takagi in \cite{T2} classified homogeneous real hypersurfaces in $\mathbb{C}P^{n}$ and Berndt in \cite{Ber} classified Hopf hypersurfaces with constant principal curvatures in $\mathbb{C}H^{n}$. Let $M$ be a real hypersurface in $M_{n}(c)$, $c\neq0$. Then we state the following theorems due to  Okumura \cite{Ok} for $\mathbb{C}P^{n}$ and Montiel and Romero \cite{MR} for $\mathbb{C}H^{n}$ respectively.\\

\begin{theorem}
Let M be a real hypersurface of $M_{n}(c)$ , $n\geq2$, $c\neq0$. If it satisfies
$A\varphi-\varphi A=0$, then M is locally congruent to one of the
following hypersurfaces:
 \begin{itemize}
    \item  In case $\mathbb{C}P^{n}$\\
    $(A_{1})$   a geodesic hypersphere of radius r , where
    $0<r<\frac{\pi}{2}$,\\
    $(A_{2})$  a tube of radius r over a totally geodesic
    $\mathbb{C}P^{k}$,$(1\leq k\leq n-2)$, where $0<r<\frac{\pi}{2}.$
    \item In case $\mathbb{C}H^{n}$\\
    $(A_{0})$   a horosphere in $ \mathbb{C}H^{n}$, i.e a Montiel tube,\\
    $(A_{1})$  a geodesic hypersphere or a tube over a hyperplane $\mathbb{C}H^{n-1}$,\\
    $(A_{2}) $  a tube over a totally geodesic $\mathbb{C}H^{k}$ $(1\leq k\leq n-2)$.
  \end{itemize}
\end{theorem}

Since 2006 many authors have studied real hypersurfaces whose structure Jacobi operator is parallel $(\nabla l=0)$. Ortega, Perez and
Santos \cite{OPS} proved the nonexistence of real hypersurfaces in non-flat complex space form with parallel structure Jacobi operator
$\nabla l=0$. Perez, Santos and Suh \cite{PSaSuh} continuing the work of \cite{OPS} considered a weaker condition ($\mathbb{D}$-parallelness),
that is $\nabla_{X}l=0$ for any vector field $X$ orthogonal to $\xi$. They proved the non-existence of such real hypersurfaces in
$\mathbb{C}P^{m}$, $m\geq3$.

Kim and Ki in \cite{KK}  classified real hypersurfaces if $\nabla_{\xi}l=0$ and $S\varphi=\varphi S$. Ki and Liu \cite{KL} proved that real
hypersurfaces satisfying $\nabla_{\xi}l=0$ and $lS=Sl$ are Hopf hypersurfaces provided that the scalar curvature is non-negative.
Ki, et.al. in \cite{KPSaSuh}  classified real hypersurfaces satisfying $\nabla_{\xi}l=0$ and $\nabla_{\xi}S=0$. Kim et.al. in \cite{KKK} studied the real hypersurfaces
satisfying $g(\nabla_{\xi}\xi,\nabla_{\xi}\xi)=\mu^{2}=$const, $6\mu^{2}+\frac{c}{4}\neq0$ and classified those whose $l$ is
$\xi-$parallel. Cho and Ki \cite{CK1} classified real hypersurfaces satisfying A$l=l$A and $\nabla_{\xi}l=0$.

Recently Ivey and Ryan, in \cite{IR} studied real hypersurfaces in $M_{2}(c)$.

Motivated by all the above conclusions we study real hypersurfacs in $\mathbb{C}P^{2}$ or
$\mathbb{C}H^{2}$ equipped with $\xi$-parallel structure Jacobi operator, i.e. $\nabla_{\xi}l=0$. More precisely, the following relation holds:
\begin{eqnarray}
(\nabla_{\xi}l)X=0.
\end{eqnarray}
We prove the following theorem

\begin{pro}
Let M be a connected real hypersurface in $\mathbb{C}P^{2}$
or $\mathbb{C}H^{2}$ with $\xi$-parallel structure Jacobi
operator. Then M is a Hopf hypersurface. Further, if
$\eta(A\xi)\neq0$, then :
\begin{itemize}
    \item in the case of  $\mathbb{C}P^{2}$, $M$ is locally congruent to\\
    a geodesic sphere, where
    $0<r<\frac{\pi}{2}$ and $r\neq\frac{\pi}{4}$,
    \item in the case of  $\mathbb{C}H^{2}$, $M$ is locally congruent \\
    to a horosphere,\\
    or to a geodesic sphere\\
    or to a tube over the hyperplane $\mathbb{C}H^{1}$.\\
\end{itemize}
\end{pro}

\section{Preliminaries}
Throughout this paper all manifolds, vector fields e.t.c. are assumed to be of class $C^{\infty}$ and all manifolds are assumed to be connected. Furthermore, the real hypersurfaces are supposed to be oriented and without boundary.

Let $M$ be a real hypersurface immersed in a nonflat complex space form $(M_{n}(c),G)$ with almost complex structure J
of constant holomorphic sectional curvature $c$. Let $N$ be a unit normal vector field on $M$ and $\xi=-JN$. For a vector field $X$ tangent to $M$ we can write $JX=\varphi (X)+\eta(X)N$, where $\varphi X$ and $\eta(X)N$ are the tangential and the normal
component of $JX$ respectively. The Riemannian connection $\overline{\nabla}$ in $M_{n}(c)$ and $\nabla$ in $M$ are related
for any vector fields $X$, $Y$ on $M$:
$$\overline{\nabla}_{Y}X=\nabla_{Y}X+g(AY,X)N$$
$$\overline{\nabla}_{X}N=-AX$$
where $g$ is the Riemannian metric on $M$ induced from G of $M_{n}(c)$ and A is the shape operator of $M$ in $M_{n}(c)$. $M$ has an almost
contact metric structure $(\varphi,\xi,\eta)$ induced from J on $M_{n}(c)$ where $\varphi$ is a (1,1) tensor field and $\eta$ a
1-form on $M$ such that (see \cite{Bl} $$g(\varphi X,Y)=G(JX,Y),\hspace{20pt}\eta(X)=g(X,\xi)=G(JX,N).$$
Then we have
\begin{eqnarray}
\varphi^{2}X=-X+\eta(X)\xi,\hspace{20pt}
\eta\circ\varphi=0,\hspace{20pt} \varphi\xi=0,\hspace{20pt}
\eta(\xi)=1
\end{eqnarray}
\begin{eqnarray}\hspace{20pt}
g(\varphi X,\varphi
Y)=g(X,Y)-\eta(X)\eta(Y),\hspace{10pt}g(X,\varphi Y)=-g(\varphi
X,Y)
\end{eqnarray}
\begin{eqnarray}
\nabla_{X}\xi=\varphi
AX,\hspace{20pt}(\nabla_{X}\varphi)Y=\eta(Y)AX-g(AX,Y)\xi
\end{eqnarray}
    Since the ambient space is of constant holomorphic sectional curvature $c$, the equations of Gauss and Codazzi for any vector
fields $X$, $Y$, $Z$ on $M$ are respectively given by
\begin{eqnarray}
R(X,Y)Z=\frac{c}{4}[g(Y,Z)X-g(X,Z)Y+g(\varphi Y ,Z)\varphi
X\end{eqnarray} $$-g(\varphi X,Z)\varphi Y-2g(\varphi X,Y)\varphi
Z]+g(AY,Z)AX-g(AX,Z)AY$$
\begin{eqnarray}
\hspace{10pt}
(\nabla_{X}A)Y-(\nabla_{Y}A)X=\frac{c}{4}[\eta(X)\varphi
Y-\eta(Y)\varphi X-2g(\varphi X,Y)\xi]
\end{eqnarray}
where $R$ denotes the Riemannian curvature tensor on $M$.

    For every point $P\epsilon M$, the tangent space $T_{P}M$ can be decomposed as following:
$$T_{P}M=span\{\xi\}\oplus ker{\eta}$$
where $ker(\eta)=\{X\;\;\epsilon\;\; T_{P}M:\eta(X)=0\}$.
 Due to the above decomposition,the vector field $A\xi$ is decomposed as follows:
 \begin{eqnarray}
 A\xi=\alpha\xi+\beta U
 \end{eqnarray}
 where $\beta=|\varphi\nabla_{\xi}\xi|$ and
 $U=-\frac{1}{\beta}\varphi\nabla_{\xi}\xi\;\epsilon\;ker(\eta)$, provided
 that $\beta\neq0$.

\section{Auxiliary relations}

Let $M$ be a real hypersurfaces in $\mathbb{C}P^{2}$ or $\mathbb{C}H^{2}$, i.e. $M_{2}(c)$, $c\neq0$. We  consider the open subset $\mathcal{N}$ of $M$ such that:
$$\mathcal{N}=\{P\;\;\epsilon\;\;M:\;\beta\neq0,\;\;\mbox{in a neighborhood of P}\}.$$
Furthermore, we consider $\mathcal{V}$, $\Omega$ open subsets of $\mathcal{N}$ such that:
$$\mathcal{V}=\{P\;\;\epsilon\;\;\mathcal{N}:\alpha=0,\;\;\mbox{in a neighborhood of P}\},$$
$$\Omega=\{P\;\;\epsilon\;\;\mathcal{N}:\alpha\neq0,\;\;\mbox{in a neighborhood of P}\},$$
where $\mathcal{V}\cup\Omega$ is open and dense in the closure of $\mathcal{N}$.

\begin{lemma}
Let M be a real hypersurface in $M_{2}(c)$, equipped with $\xi$-parallel structure Jacobi operator. Then
$\mathcal{V}$ is empty.
\end{lemma}
\textbf{Proof:} Let $\{U,\varphi U,\xi\}$ be a local orthonormal basis on $\mathcal{V}$. The
relation (2.6) takes the form $A\xi=\beta U$. The first relation of (2.3) for $X=\xi$, taking into account the latter implies
\begin{eqnarray}
\nabla_{\xi}\xi=\beta\varphi U.\nonumber\
\end{eqnarray}
Relation (1.1) for $X=\xi$, because of the above relation yields:
\begin{eqnarray}
\nabla_{\xi}(l\xi)=l\nabla_{\xi}\xi\Rightarrow \beta\varphi U=0,\nonumber\
\end{eqnarray}
which leads to a contradiction and this completes the proof of Lemma 3.1.
\qed
\\

In what follows we work on $\Omega$, where $\alpha\neq0$ and $\beta\neq0$.

\begin{lemma}
Let M be a real hypersurface in $M_{2}(c)$, equipped with $\xi$-parallel structure Jacobi operator. Then the following relations hold in $\Omega$:
\begin{eqnarray}
\hspace{-140pt}AU=(\frac{\beta^{2}}{\alpha}-\frac{c}{4\alpha}+\frac{\kappa}{\alpha})U+\beta\xi,\;\;\;\;
A\varphi U=-\frac{c}{4\alpha}\varphi U
\end{eqnarray}
\begin{eqnarray}
\hspace{-100pt}\nabla_{\xi}\xi=\beta\varphi U,\;\;\;
\nabla_{U}\xi=(\frac{\beta^{2}}{\alpha}-\frac{c}{4\alpha}+\frac{\kappa}{\alpha})\varphi
U,\;\;\; \nabla_{\varphi U}\xi=\frac{c}{4\alpha}U
\end{eqnarray}
\begin{eqnarray}
\hspace{-100pt}\nabla_{\xi}U=\kappa_{1}\varphi U,\;\;\;
\nabla_{U}U=\kappa_{2}\varphi U,\;\;\; \nabla_{\varphi
U}U=\kappa_{3}\varphi U-\frac{c}{4\alpha}\xi
\end{eqnarray}
\begin{equation}
\nabla_{\xi}\varphi U=-\kappa_{1}U-\beta\xi,\;\;\;
\nabla_{U}\varphi
U=-\kappa_{2}U-(\frac{\beta^{2}}{\alpha}-\frac{c}{4\alpha}+\frac{\kappa}{\alpha})\xi,\;\;\;
\nabla_{\varphi U}\varphi U=-\kappa_{3}U
\end{equation}
\begin{eqnarray}
\hspace{-180pt}\kappa\kappa_{1}=0,\;\;\;   (\xi\kappa)=0,
\end{eqnarray}
where $\kappa,\kappa_{1},\kappa_{2},\kappa_{3}$ are smooth
functions on M.
\end{lemma}
\textbf{Proof:} Let $\{U,\varphi U,\xi\}$ be a local orthonormal basis of $\Omega$.

The first relation of (2.3) for $X=\xi$ implies: $\nabla_{\xi}\xi=\beta\varphi U$ and so relation (1.1) for $X=\xi$, taking into account the latter, gives:
\begin{eqnarray}
l\varphi U=0.
\end{eqnarray}
Relation (2.4) for $X=\varphi U$ and $Y=Z=\xi$ gives: $l\varphi U=\frac{c}{4}\varphi U+\alpha A\varphi U$, which because of (3.6) implies the second of (3.1). Relation (2.4) for $X=U$ and $Y=Z=\xi$, we have:
\begin{eqnarray}
lU=\frac{c}{4}U+\alpha AU-\beta A\xi\
\end{eqnarray}
The scalar products of (3.7) with $\varphi U$ and $U$, because of (2.6) and the second of (3.1) imply the first of (3.1), where $\kappa=g(lU,U)$.

The first relation of (2.3), for $X=U$ and $X=\varphi U$. taking into consideration relations (3.1), gives the rest of relation (3.2).

From the well known relation: $Xg(Y,Z)=g(\nabla_{X}Y,Z)+g(Y,\nabla_{X}Z)$ for $X,Y,Z\;\;
\epsilon$
$\{\xi,U,\varphi U\}$ we obtain (3.3) and (3.4), where $\kappa_{1},\kappa_{2},\kappa_{3}$ are
smooth functions in $\Omega$.

On the other hand
 \begin{eqnarray}
&&\xi\kappa=\xi g(lU,U)\nonumber\\
&&\Rightarrow \xi\kappa=g(\nabla_{\xi}(lU),U)+g(lU,\nabla_{\xi}U)\nonumber\\
&&\Rightarrow \xi\kappa=g((\nabla_{\xi}l)U+l(\nabla_{\xi}U),U)+g(lU,\nabla_{\xi}U)\nonumber\\
&&\Rightarrow \xi\kappa=g(l(\nabla_{\xi}U),U)+g(lU,\nabla_{\xi}U)\nonumber\
\end{eqnarray}
The above relation because of (3.3), (3.6) and (3.7) yields:
\begin{eqnarray}
\xi\kappa=g(\kappa_{1}l\varphi U,U)+g(lU,\kappa_{1}\varphi
U)\Rightarrow \xi\kappa=0\nonumber\
\end{eqnarray}
On the other hand:
\begin{eqnarray}
&&\xi g(l\varphi U,U)=0\nonumber\\
&&\Rightarrow g(\nabla_{\xi}(l\varphi U),U)+g(l\varphi U,\nabla_{\xi}U)=0\nonumber\\
&&\Rightarrow g((\nabla_{\xi}l)\varphi U+l(\nabla_{\xi}\varphi U),U)+g(l\varphi U,\nabla_{\xi}U)=0\nonumber\
\end{eqnarray}
From the above equation because of (1.1), (2.6), (3.4), (3.6) and $\kappa=g(lU,U)$ we obtain:
\begin{eqnarray}
&&g(l(-\kappa_{1}U-\beta\xi),U)=0\Rightarrow \kappa\kappa_{1}=0\nonumber\
\end{eqnarray}
\qed

Relation (2.5) for $X$ $\epsilon$  $\{U,\varphi U\}$ and $Y=\xi$, because of Lemma 3.2 yields:
\begin{eqnarray}
U\beta&=&\xi(\frac{\beta^{2}}{\alpha}-\frac{c}{4\alpha}+\frac{\kappa}{\alpha})\\
U\alpha&=&\xi\beta\\
\frac{\beta^{2}\kappa_{1}}{\alpha}&=&\kappa+\beta\kappa_{2}+\frac{c}{4\alpha}(\frac{\kappa}{\alpha}-\frac{c}{4\alpha}+\frac{\beta^{2}}{\alpha})\\
(\varphi U)\beta&=&\frac{\kappa_{1}\beta^{2}}{\alpha}+\beta^{2}+\frac{c}{4\alpha}(\frac{\beta^{2}}{\alpha}-\frac{c}{4\alpha}+\frac{\kappa}{\alpha})\\
\xi\alpha&=&\frac{4\alpha^{2}\kappa_{3}\beta}{c}\\
(\varphi U)\alpha&=&\beta(\kappa_{1}+\alpha+\frac{3c}{4\alpha})
\end{eqnarray}
Furthermore, relation (2.5), for $X=U$ and $Y=\varphi U$, due to Lemma 3.2 and (3.10), implies:
\begin{eqnarray}
(\varphi U)\kappa&=&-\frac{c\beta\kappa_{1}}{4\alpha}+\kappa\beta+\kappa\kappa_{2}-c\beta\\
U\alpha&=&\frac{4\kappa_{3}\alpha}{c}(\beta^{2}+\kappa)
\end{eqnarray}
Using the relations (3.9)-(3.15) and Lemma 3.2 we obtain:
 \begin{eqnarray}
&&[U,\xi](\frac{c}{4\alpha})=(\nabla_{U}\xi-\nabla_{\xi}U)\frac{c}{4\alpha}\nonumber\\
&&\Rightarrow [U,\xi](\frac{c}{4\alpha}) =-\frac{c\beta}{4\alpha^{2}}(\frac{\beta^{2}}{\alpha}-\frac{c}{4\alpha}+\frac{\kappa}{\alpha}\frac{}{}-\kappa_{1})(\kappa_{1}+\alpha+\frac{3c}{4\alpha})
\end{eqnarray}
\begin{eqnarray}
&& [U,\xi](\frac{c}{4\alpha})=(U(\xi\frac{c}{4\alpha}))-(\xi(U\frac{c}{4\alpha}))\nonumber\\
&&\Rightarrow [U,\xi](\frac{c}{4\alpha})=-\beta\kappa_{3}^{2}-\beta(U\kappa_{3})+(\frac{\beta^{2}}{\alpha}+\frac{\kappa}{\alpha})(\xi\kappa_{3})
\end{eqnarray}
Similarly:
\begin{eqnarray}
&&[U,\varphi
U](\frac{c}{4\alpha})=\kappa_{2}\kappa_{3}(\frac{\beta^{2}}{\alpha}+\frac{\kappa}{\alpha})+\beta\kappa_{3}(\frac{\kappa}{\alpha}-\frac{c}{2\alpha}+\frac{\beta^{2}}{\alpha})\nonumber\\
&&+\frac{c\beta\kappa_{3}}{4\alpha^{2}}(\kappa_{1}+\alpha+\frac{3c}{4\alpha})
\end{eqnarray}
\begin{eqnarray}
&&[U,\varphi
U](\frac{c}{4\alpha})=\frac{2\kappa_{3}\beta^{3}\kappa_{1}}{\alpha^{2}}+\frac{\kappa_{3}\beta^{3}}{\alpha}+\frac{5c\kappa_{3}\beta}{4\alpha^{3}}(\beta^{2}+\kappa)-\frac{c\beta\kappa_{1}\kappa_{3}}{4\alpha^{2}}\nonumber\\
&&-\frac{c\beta\kappa_{3}}{4\alpha}-\frac{5c^{2}\beta\kappa_{3}}{16\alpha^{3}}-\frac{c\beta}{4\alpha^{2}}(U\kappa_{1})-\frac{\beta\kappa\kappa_{3}}{\alpha}+\frac{\kappa_{3}}{\alpha}((\varphi U)\kappa)\nonumber\\
&&+(\frac{\beta^{2}}{\alpha}+\frac{\kappa}{\alpha})((\varphi U)\kappa_{3})
\end{eqnarray}
\begin{eqnarray}
[\varphi
U,\xi](\frac{c}{4\alpha})=-\kappa_{3}(\kappa_{1}+\frac{c}{4\alpha})(\frac{\beta^{2}}{\alpha}+\frac{\kappa}{\alpha})-\beta^{2}\kappa_{3}
\end{eqnarray}
\begin{eqnarray}
&&[\varphi
U,\xi](\frac{c}{4\alpha})=-\frac{2\kappa_{1}\kappa_{3}\beta^{2}}{\alpha}-\kappa_{3}\beta^{2}-\frac{7c\beta^{2}\kappa_{3}}{4\alpha^{2}}+\frac{c^{2}\kappa_{3}}{16\alpha^{2}}+\frac{c\kappa\kappa_{3}}{2\alpha^{2}}\nonumber\\
&&-\beta(\varphi U)\kappa_{3}+\kappa\kappa_{3}+\frac{c\beta}{4\alpha^{2}}(\xi\kappa_{1}).
\end{eqnarray}

Due to the first relation of (3.5), we consider $\Omega_{1}$ the open subset of $\Omega$ such that:
$$\Omega_{1}=\{P\;\;\epsilon\;\;\Omega:\kappa_{1}\neq0,\;\;in\;\;a\;\;neighborhood\;\;of\;\;P\}.$$
So in $\Omega_{1}$, we have: $\kappa=0$.

In $\Omega_{1}$ relation (3.14), since $\kappa=0$, yields:
\begin{eqnarray}
\kappa_{1}=-4\alpha
\end{eqnarray}
and from relation (3.10), taking into account (3.22), we get:
\begin{eqnarray}
\kappa_{2}=-4\beta-\frac{c\beta}{4\alpha^{2}}+\frac{c^{2}}{16\alpha^{2}\beta}
\end{eqnarray}
From (3.20) and (3.21), using (3.12), (3.22) and (3.23) we obtain:
\begin{eqnarray}
\beta(\varphi U)\kappa_{3}=-\frac{3c\beta^{2}\kappa_{3}}{2\alpha^{2}}+\frac{c^{2}\kappa_{3}}{16\alpha^{2}}
\end{eqnarray}
From (3.18), (3.19), using (3.15), (3.22), (3.23) and (3.24), we obtain:
\begin{eqnarray}
\kappa_{3}(4\alpha^{2}-c)=0.
\end{eqnarray}
Because of (3.25), let $\Omega'_{1}$ be the open subset of $\Omega_{1}$ such that:
$$\Omega'_{1}=\{P\;\;\epsilon\;\;\Omega_{1}:\kappa_{3}\neq0,\;\;in\;\;a\;\;neighborhood\;\;of\;\;P\}.$$
So in $\Omega'_{1}$ we obtain: $c=4\alpha^{2}$. Differentiation of the latter with respect to
$\xi$, implies $\xi\alpha=0$ which because of (3.12) leads to $\kappa_{3}=0$, which is impossible. So $\Omega'_{1}$ is empty and $\kappa_{3}=0$ in $\Omega_{1}$.

\begin{lemma}
Let M be a real hypersurface in $M_{2}(c)$, equipped with $\xi$-parallel structure Jacobi operator. Then $\Omega_{1}$ is empty.
\end{lemma}
\textbf{Proof:} We resume that in $\Omega_{1}$ we have:
\begin{eqnarray}
\kappa=\kappa_{3}=0
\end{eqnarray}
and relations (3.22), (3.23) and (3.24) hold.

Relations (3.8), (3.9), (3.12) and (3.15), because of (3.5) and (3.26), yield:
\begin{eqnarray}
U\alpha=U\beta=\xi\alpha=\xi\beta=0
\end{eqnarray}
In $\Omega_{1}$, combining (3.16) and (3.17) and taking into account (3.22) and (3.26), we obtain:
\begin{eqnarray}
(\frac{c}{4\alpha}-\alpha)(\frac{\beta^{2}}{\alpha}-\frac{c}{4\alpha}+4\alpha)=0
\end{eqnarray}
Owing to (3.28), let $\Omega_{11}$ be the open subset of $\Omega_{1}$, such that:
$$\Omega_{11}=\{P\;\;\epsilon\;\;\Omega_{1}:c\neq4\alpha^{2},\;\;\mbox{in a neighborhood of P}\}.$$
From (3.28) in $\Omega_{11}$, we have: $4\alpha=-\frac{\beta^{2}}{\alpha}+\frac{c}{4\alpha}$. Differentiation of the latter along $\varphi U$, because of (3.11), (3.13), (3.22), (3.26) and the last relation yields $c=0$, which is impossible. Hence, $\Omega_{11}$ is empty.

So in $\Omega_{1}$ the relation $c=4\alpha^{2}$ holds. Due to the last relation and (3.22), the relation (3.11) becomes:
\begin{eqnarray}
(\varphi U)\beta=-(\alpha^{2}+2\beta^{2}).
\end{eqnarray}
From (3.27) we have $[U,\xi]\beta=U(\xi\beta)-\xi(U\beta)\Rightarrow [U,\xi]\beta=0$. On the other hand, from (3.2) and (3.3) we obtain
 $[U,\xi]\beta=(\nabla_{U}\xi-\nabla_{\xi}U)\beta\Rightarrow [U,\xi]\beta=\frac{1}{\alpha}(3\alpha^{2}+\beta^{2})(\varphi U)\beta$. The
 last two relations imply $(\varphi U)\beta=0$. Therefore, from (3.29) we obtain $\alpha^{2}+2\beta^{2}=0$, which is a contradiction. Hence, $\Omega_{1}$ is empty.
\qed
\\

Since $\Omega_{1}$ is empty, in $\Omega$ we have $\kappa_{1}=0$. So from relations (3.20) and (3.21) we obtain:
$$\beta(\varphi
U)\kappa_{3}=\frac{\kappa_{3}}{16\alpha^{2}}[c^{2}-24c\beta^{2}+12c\kappa+16\alpha^{2}\kappa].$$
Furthermore, the combination of relations (3.18) and (3.19), using (3.10) and (3.14), implies:
$$(\beta^{2}+\kappa)(\varphi
U)\kappa_{3}=\frac{c\beta\kappa_{3}}{16\alpha^{2}}[16\alpha^{2}-24(\beta^{2}+\kappa)+9c].$$
From the last two relations we obtain:
\begin{eqnarray}
\kappa_{3}[c^{2}\kappa+12c\kappa^{2}+12c\beta^{2}\kappa+16\alpha^{2}\beta^{2}\kappa-16c\alpha^{2}\beta^{2}+16\alpha^{2}\kappa^{2}-8c^{2}\beta^{2}]=0\nonumber\
\end{eqnarray}
Due to the above relation, we consider $\Omega_{2}$ the open subset of $\Omega$, such that:
$$\Omega_{2}=\{P\;\;\epsilon\;\;\Omega:\kappa_{3}\neq0,\;\;\mbox{in a neighborhood of P}\},$$
so in $\Omega_{2}$ the following relation holds:
\begin{eqnarray}
c^{2}\kappa+12c\kappa^{2}+12c\beta^{2}\kappa+16\alpha^{2}\beta^{2}\kappa-16c\alpha^{2}\beta^{2}+16\alpha^{2}\kappa^{2}-8c^{2}\beta^{2}=0.
\end{eqnarray}
Differentiating (3.30) with respect to $\xi$ and using (3.5), (3.9), (3.12) and (3.15) we obtain:
\begin{eqnarray}
&&8\alpha^{2}\beta^{2}\kappa-8c\alpha^{2}\beta^{2}+8\alpha^{2}\kappa^{2}+3c\beta^{2}\kappa-2c^{2}\beta^{2}+3c\kappa^{2}-4c\alpha^{2}\kappa\nonumber\\
&&-2c^{2}\kappa=0
\end{eqnarray}
From (3.30) and (3.31) we obtain:
\begin{eqnarray}
5c\kappa+6\kappa^{2}+6\beta^{2}\kappa-4c\beta^{2}+8\alpha^{2}\kappa=0.
\end{eqnarray}
Differentiating (3.32) with respect to $\xi$ and using (3.5), (3.9), (3.12) and (3.15) we have:
$4\kappa\alpha^{2}=(2c-3\kappa)(\beta^{2}+\kappa).$ The last relation with (3.32) imply: $\kappa=0$. Substituting the latter in (3.30) gives $c=-2\alpha^{2}$. Differentiation of the last relation with respect to $\varphi U$ and taking into account (3.13), $c=-2\alpha^{2}$ and $\kappa_{1}=0$ results in $\alpha=0$, which is impossible.

So $\Omega_{2}$ is empty and in $\Omega$ we get: $\kappa_{3}=0$.

\begin{lemma}
Let M be a real hypersurface in $M_{2}(c)$, equipped with $\xi$-parallel structure Jacobi operator. Then $\Omega$ is empty.
\end{lemma}
\textbf{Proof:} We resume that in $\Omega$ the following relation holds:
\begin{eqnarray}
\kappa_{1}=\kappa_{3}=0.
\end{eqnarray}
Relations (3.8), (3.9), (3.12), (3.15), because of (3.5) and (3.33) yield:
\begin{eqnarray}
U\alpha=U\beta=\xi\alpha=\xi\beta=0
\end{eqnarray}
In $\Omega$ the combination of (3.16), (3.17) and taking into account (3.33), implies:
\begin{eqnarray}
(4\alpha^{2}+3c)(\beta^{2}+\kappa-\frac{c}{4})=0
\end{eqnarray}
Due to (3.35), we consider $\Omega_{3}$ the open subset of $\Omega$ such that:
$$\Omega_{3}=\{P\;\;\epsilon\;\;\Omega:\beta^{2}+\kappa\neq\frac{c}{4},\;\;\mbox{in a neighborhood of P}\}.$$
So in $\Omega_{3}$ the following relation holds:
\begin{eqnarray}
c=-\frac{4\alpha^{2}}{3}.
\end{eqnarray}
Differentiation of (3.36) with respect to $\varphi U$ implies:
\begin{eqnarray}
(\varphi U)\alpha=0.
\end{eqnarray}
Because of (3.34) we have $[U,\xi]\beta=U(\xi\beta)-\xi(U\beta)\Rightarrow [U,\xi]\beta=0$. On the other hand  due to
(3.2), (3.3) and (3.33) we get $[U,\xi]\beta=(\nabla_{U}\xi-\nabla_{\xi}U)\beta\Rightarrow [U,\xi]\beta=(\frac{\beta^{2}}{\alpha}-\frac{c}{4\alpha}+\frac{\kappa}{\alpha})(\varphi
U)\beta$. Combination of the last relations imply:
\begin{eqnarray}
(\varphi U)\beta=0
\end{eqnarray}
From (3.11), owing to (3.33), (3.36) and (3.38) yields $2\beta^{2}=\kappa+\frac{\alpha^{2}}{3}$. Differentiation of the
last relation with respect to $\varphi U$ and taking into account (3.37) and (3.38) imply
$(\varphi U)\kappa=0$. So from (3.14), because of the latter and (3.33), we obtain $\kappa(\beta+\kappa_{2})=c\beta$. The combination of the latter with (3.10) and taking into account (3.33), (3.36) and $2\beta^{2}=\kappa+\frac{a^{2}}{3}$ imply:
\begin{eqnarray}
\alpha^{2}=18\beta^{2}\hspace{20pt}\kappa_{2}=5\beta\hspace{20pt}\kappa=-4\beta^{2}
\end{eqnarray}
The relations of Lemma 3.2 in $\Omega_{3}$, because of (3.36) and (3.39) become:
\begin{eqnarray}
&&AU=\frac{\alpha}{6}U+\beta\xi,\;\;\;A\varphi U=\frac{\alpha}{3}\varphi U\\
&&\nabla_{\xi}\xi=\beta\varphi U,\;\;\;\nabla_{U}\xi=\frac{\alpha}{6}\varphi U,\;\;\;\nabla_{\varphi U}\xi=-\frac{\alpha}{3}U,\\
&&\nabla_{\xi}U=0,\;\;\;\nabla_{U}U=5\beta\varphi U,\;\;\;\nabla_{\varphi U}U=\frac{\alpha}{3}\xi,\\
&&\nabla_{\xi}\varphi U=-\beta\xi,\;\;\;\nabla_{U}\varphi U=-5\beta U-\frac{\alpha}{6}\xi,\;\;\;\nabla_{\varphi U}\varphi U=0.
\end{eqnarray}

The relation (2.4), because of (3.36), (3.39) and (3.40) implies: $R(U,\varphi U)U=23\beta^{2}\varphi U$. On the other hand
$R(X,Y)Z=\nabla_{X}\nabla_{Y}Z-\nabla_{Y}\nabla_{X}Z-\nabla_{[X,Y]}Z$, because of (3.34), (3.36), (3.39) and(3.41)-(3.43) yields: $R(U,\varphi
U)U=26\beta^{2}\varphi U$. The combination the last two relations implies $\beta=0$, which is impossible in $\Omega_{3}$.

So $\Omega_{3}$ is empty and in $\Omega$ the following relation holds
\begin{eqnarray}
\beta^{2}+\kappa=\frac{c}{4}.
\end{eqnarray}
In $\Omega$ (3.10) becomes:
\begin{eqnarray}
\kappa+\beta\kappa_{2}=0.
\end{eqnarray}

Differentiating (3.44) with respect to $\varphi U$ and using (3.11), (3.14), (3.33), (3.44) and (3.45) we obtain:
$\beta^{2}=-\frac{c}{4}$. Differentiation of the last relation along $\varphi U$ implies $(\varphi U)\beta=0$, which because of (3.11), (3.33) and (3.44) yields $\beta=0$, which is a contradiction. Therefore, $\Omega$ is empty and this completes the proof of Lemma 3.4.
\qed
\\

From Lemmas 3.1 and 3.4, we conclude that $\mathcal{N}$ is empty and we lead to the following result:
\begin{proposition}
Every real hypersurface in $M_{2}(c)$, equipped with $\xi$-parallel structure Jacobi operator, is a Hopf hypersurface.
\end{proposition}

\section{\hspace{-15pt}.\hspace{10pt}Proof of Main Theorem}
Since $M$ is a Hopf hypersurface, due to Theorem 2.1 \cite{NR1}  we have that $\alpha$ is a constant. We suppose that $\alpha\neq0$. We consider a unit vector field $e$ $\epsilon$ $\mathbb{D}$, such that $Ae=\lambda e$, then $A\varphi e=\nu\varphi e$ at some point $P$ $\epsilon$ $M$, where $\{ e, \varphi e, \xi\}$ is a local orthonormal basis. Then the following relation holds on $M$, (Corollary 2.3 \cite{NR1}):
\begin{eqnarray}
\lambda\nu=\frac{\alpha}{2}(\lambda+\nu)+\frac{c}{4}.
\end{eqnarray}
The first relation of (2.3) for $X=e$ implies:
\begin{eqnarray}
\nabla_{e}\xi=\lambda\varphi e.
\end{eqnarray}
Relation (2.4) for $X=e$ and $Y=Z=\xi$ yields:
\begin{eqnarray}
le=\frac{c}{4}e+\alpha Ae.
\end{eqnarray}
From relation (1.1) for $X=e$, we obtain:
\begin{eqnarray}
\nabla_{\xi}(le)=l\nabla_{\xi}e.
\end{eqnarray}
From (2.4) for $X=\nabla_{\xi}e$ and $Y=Z=\xi$, we get:
\begin{eqnarray}
l\nabla_{\xi}e=\frac{c}{4}\nabla_{\xi}e+\alpha A(\nabla_{\xi}e).
\end{eqnarray}
Substitution in (4.4) of (4.3) and (4.5) yields:
\begin{eqnarray}
(\nabla_{\xi}A)e=0.
\end{eqnarray}
The relation (2.5) for $X=\xi$ and $Y=e$, taking into account (4.6), we get:
\begin{eqnarray}
(\nabla_{e}A)\xi=-\frac{c}{4}\varphi e
\end{eqnarray}
Finally, the scalar product of (4.7) with $\varphi e$, taking into consideration (4.1), (4.2) and $A\varphi e=\nu\varphi e$ yields:
\begin{eqnarray}
&&g(\nabla_{e}(A\xi)-A\nabla_{e}\xi, \varphi e)=-\frac{c}{4}\nonumber\\
&& \Rightarrow \alpha\lambda=-\frac{c}{4}+\lambda\nu\Rightarrow \lambda=\nu.\nonumber\
\end{eqnarray}
Then $Ae=\lambda e$ and $A\varphi e=\lambda\varphi e$, therefore we obtain:  $$(A\varphi-\varphi A)X=0,\;\;\forall\;\;X\;\;\epsilon\;\;TM.$$
From the above relation Theorem 1.1 holds. Since $\alpha\neq0$ we can not have the geodesic sphere of radius $r=\frac{\pi}{4}$ and this completes the Proof of Main Theorem.

\section*{Acknowledgements}
The first author is granted by the Foundation Alexandros S. Onasis. Grant Nr: G ZF 044/2009-2010.

\end{document}